\documentclass[preprint,authoryear,12pt]{article}
\usepackage{graphicx,amsmath,amsfonts,amssymb,amsthm,color,slashbox,amsbsy,ragged2e,natbib,hyperref}

\newtheorem{pretheorem}{{\bf Theorem}}[section]

\newtheorem{preresult}[pretheorem]{{\bf Result }}

\newtheorem{prelemma}[pretheorem]{{\bf Lemma}}

\newtheorem{preexample}[pretheorem]{{\bf Example}}

\newtheorem{predefn}[pretheorem]{{\bf Definition  }}

\newtheorem{prenote-}[pretheorem]{{\bf Remark }}

\newenvironment{note-}{\begin{prenote-}{\hspace{-0.5em} {\bf }}}
                         {\end{prenote-}}
\newcommand{\bq}{\begin{equation}}
\newcommand{\eq}{\end{equation}}
\newcommand{\bqs}{\begin{equation*}}
\newcommand{\eqs}{\end{equation*}}
\newcommand{\bqas}{\begin{eqnarray*}}
\newcommand{\eqas}{\end{eqnarray*}}
\newcommand{\bqa}{\begin{eqnarray}}
\newcommand{\eqa}{\end{eqnarray}}

\def\bp{\boldsymbol{p}}
\def\bt{\boldsymbol{t}}
\begin{document}
\title{Goodness of fit test under progressive Type-I interval censoring}
\author{H. Nadeb$^1$,  H. Torabi$^1$, G.G. Hamedani$^2$ \\
{\small $^1$Department of Statistics, Yazd University, 89175-741, Yazd, IRAN }\\
{\small $^2$Department of Mathematics, Statistics and Computer Science, }\\
{\small Marquette University, Milwaukee, WI, USA}
}

\date{}
\maketitle

\begin{abstract}
In this paper, we propose several statistics for testing uniformity under progressive Type-I interval censoring. We obtain the critical points of these statistics and  study the power of the proposed tests  against a representative set of alternatives via simulation. Finally, we generalize our methods for continuous and completely specified distributions. 
\end{abstract}
\noindent{\bf Keywords:} Empirical reliability function; Goodness of fit test; Progressive Type-I interval censoring; Theorical reliability function.
\section{Introduction}
 Aggarwala \citeyearpar{Aggarwala} introduced Type-I interval and progressive censoring and developed the statistical inference for the
exponential distribution based on progressively Type-I interval censored data. Ng and Wang \citeyearpar{Ng} introduced the concept of progressive Type-I interval censoring to the
Weibull distribution and compared many different estimation methods for two parameters in the Weibull distribution via simulation. In general, for progressive Type-I integral censoring, relatively little work has been done.

Suppose that $n$ items are placed on a life testing problem simultaneously at time $ t_0=0$  under inspection at $m$ pre-specified times
$t_1 < t_2 <\ldots < t_m $ where $t_m$ is the scheduled time to terminate the experiment. At the $i$th inspection time, $t_i$, the number, $X_i$,
of failures within $\left(  t_{i-1},t_i\right] $ is recorded and $R_i$ surviving items are randomly removed from the life testing, for $ i=1,\ldots ,m $.
It is obvious that the number of surviving items at the time $t_i$ is $Y_i=n-\sum_{j=1}^{i} X_j-\sum_{j=1}^{i-1} R_j$. Since $Y_i$
 is a random variable and the exact number of items withdrawn should not be greater
than $Y_i$ at time schedule $t_i$, $R_i$ could be determined by the pre-specified percentage of the remaining surviving units at $t_i$ for $i=1,2,\ldots,m$. Also,  given pre-specified percentage values, $ p_1,\ldots, p_{m-1}$ and $p_m=1$, for withdrawing at
$t_1 < t_2 < \ldots < t_m$, respectively, $ R_i=\lfloor p_iy_i\rfloor$  at each inspection time $t_i$ where $ i = 1,2,\ldots,m$. Therefore, a progressively
Type-I interval censored sample can be denoted by $(X_i, R_i, t_i)$, $i = 1, 2,\ldots m$, where sample size is $ n = \mathop{\sum}\limits_{i=1}^{m} (X_i+R_i)$. Note that if
$R_i=0$, $i=1, 2,\ldots,m-1$, then the progressively Type-I interval censored sample is a Type-I interval censored sample.

Let $(X_i, R_i, t_i)$, $i=1, \ldots , m$, be progressively Type-I interval censored sample with pre-specified vector $\bp=(p_1, \ldots, p_{m-1}, 1)$ and $t_m <1$ from an unknown distribution function $F(.)$. We are interested in the hypothesis testing
\begin{equation}\label{hypothesis}
\left\{\begin{array}{l}H_0:~F(x)=x \\ H_1:~F(x)\neq x.\end{array}\right.
\end{equation} 

Most of the goodness of fit tests are based on the distance between empirical  reliability function and theoretical reliability function over the interval (0, 1). Based on progressively Type-I interval censored sample, in view of Balakrishnan et al. \citeyearpar{Bala 2010} or Balakrishnan and Cramer \citeyearpar{Bala 2015}, the reliability at $t_i$ can be estimated nonparametrically by
\bq \label{reliability}
\hat{\bar{F}}(t_i)=\prod_{j=1}^i\left(1-\frac{X_j}{\alpha_{j-1}^+}\right),
\eq
where
\bqs
\alpha_j^+=n-X_{\bullet j}-R_{\bullet j},
\eqs
and
\bqas
X_{\bullet j}&=&\sum_{k=1}^j X_k,\\
R_{\bullet j}&=&\sum_{k=1}^j R_k,\\
\eqas
which will be used to establish statistics for \eqref{hypothesis}. 

The paper is organized as follows:\\
In Section 2, we propose several statistics for testing uniformity under progressive Type-I interval censoring. In Section 3, we obtain the critical points of these statistics and then study power of the proposed tests  against a representative set of alternatives using simulation. In Section 4, we generalize these methods for continuous and completely specified distributions.

\section{Proposed tests}
Let $(X_i, R_i, t_i)$, $i=1, \ldots , m$, be progressively Type-I interval censored sample with pre-specified vector $\bp=(p_1, \ldots, p_{m-1}, 1)$ and $t_m <1$. It is clear that, under $H_0$, we have:
\bqs
\bar{F}(t_i)=1-t_i.
\eqs
Now we consider the difference between empirical  reliability function and theoretical reliability function and define:
\bqs
D_i=\hat{\bar{F}}(t_i)-\bar{F}(t_i).
\eqs
Currently, based on $D_i$, we introduce the goodness of fit test statistics as follows:
\bqas
C^+&=&\mathop {\max}\limits_{1\leq i \leq m}(D_i),\qquad C^-=\mathop {\max} \limits_{1\leq i \leq m}(-D_i),\qquad C=\max (C^+, C^-),\\
 K&=&C^++C^-,\qquad T^{(1)}=\frac{1}{m}\sum_{i=1}^m D_i^2,\qquad ~T^{(2)}=\frac{1}{m}\sum_{i=1}^m |D_i|.
\eqas 

If the null hypothesis is true, we expect the deviation $D_i$ to be small and consequently the above test statistics to be small. Hence, we may reject the null hypothesis if the above test statistics exceed the corresponding upper-tail null critical values.

It should be noted that Pakyari and Balakrishnan \citeyearpar{Pakyari} used statistics similar to the above statistics for goodness of fit test under Type-II progressive censoring, with this difference that they used the difference between the observed value and the expectation of $i$th Type-II progressively order statistic from uniform(0,1) distribution. 

\section{Simulation  study}
In this section, we assess the power of the proposed tests by comparing the simulated power values. We generated 20,000 random samples for different choices of sample sizes and progressive censoring schemes for determining the power. For comparative purposes, we consider two vectors for inspection time as follows: 
\bqas
\bt_1&=&(0,0.1,0.2,0.3,0.4,0.5),\\
\bt_2&=&(0,0.05,0.1,0.2,0.45,0.5),\\
\eqas
and two percentage vectors as follows:
\bqas
\bp_1&=&(0.25,0.25,0.5,0.5,1),\\
\bp_2&=&(0.5,0.5,0.25,0.25,1).\\
\eqas
So, we consider three families of alternative distributions with support in [0,1]. They are defined by the following CDFs:

(a) Lehmann alternatives,
\bqs
F_{\alpha}(x)=x^{\alpha},\quad 0\leq x \leq 1,~\alpha>0;
\eqs

(b) centered distributions having a U-shaped PDF, for $\beta \in (0,1)$ and wedge-shaped PDF, for $\beta >1$, 
\bqs
F_{\beta}(x)=\left\{\begin{array}{l}\frac{1}{2}(2x)^{\beta}\qquad \qquad \qquad 0\leq x \leq \frac{1}{2}, \\ 1-\frac{1}{2}\{2(1-x)\}^{\beta} \quad~ \frac{1}{2}\leq x \leq 1;\end{array}\right.
\eqs

(c) compressed uniform alternatives,
\bqs
F_{\gamma}(x)=\frac{x-\gamma}{1-2\gamma}, \quad \gamma \leq x \leq 1-\gamma,
\eqs
where $0\leq \gamma <\frac{1}{2}$. See Fortiana and Grane \citeyearpar{Fortiana}.

As an illustration of the tests we depict only the power functions at $0.05$ significance level for $n=40$ (because in progressive Type-I interval censoring problems the sample size is relatively large) for every censoring scheme. We take the critical regions computed and listed in Table \ref{t1} via simulation. The points computed for each power curve are estimated by the relative frequency of every statistic in the critical region for 20,000 simulated samples of the alternative distribution under progressive Type-I interval censoring.

\begin{table}[ht!]
\begin{center}
\begin{tabular}{c c cccccc}
\hline $\bt$ &$\bp$& $C^+$ & $C^-$&$C$&$K$&$T^{(1)}$ &$T^{(2)}$ \\
\hline
\hline
$\bt_1$&$\bp_1$&$0.2361$&$0.2597$&$0.3140$&$0.3157$&$0.0284$&$0.1420$\\
 \hline
 $\bt_1$&$\bp_2$&$0.2775$&$0.3111$&$0.3412$&$0.3513$&$0.0409$&$0.1700$\\
\hline
$\bt_2$&$\bp_1$&$0.2470$&$0.2417$&$0.3066$&$0.3214$&$0.0310$&$0.1378$\\
\hline
$\bt_2$&$\bp_2$&$0.3000$&$0.3132$&$0.3500$&$0.3750$&$0.0440$&$0.1655$\\
\hline
\end{tabular}
\caption{\small{Simulated critical values of $C^+$, $C^-$, $C$, $K$, $T^{(1)}$ and $T^{(2)}$ at $0.05$ significance level.}}\label{t1}
\end{center}
\end{table}

\begin{figure}[!th] 
\includegraphics[width=4.5cm,height=5cm]{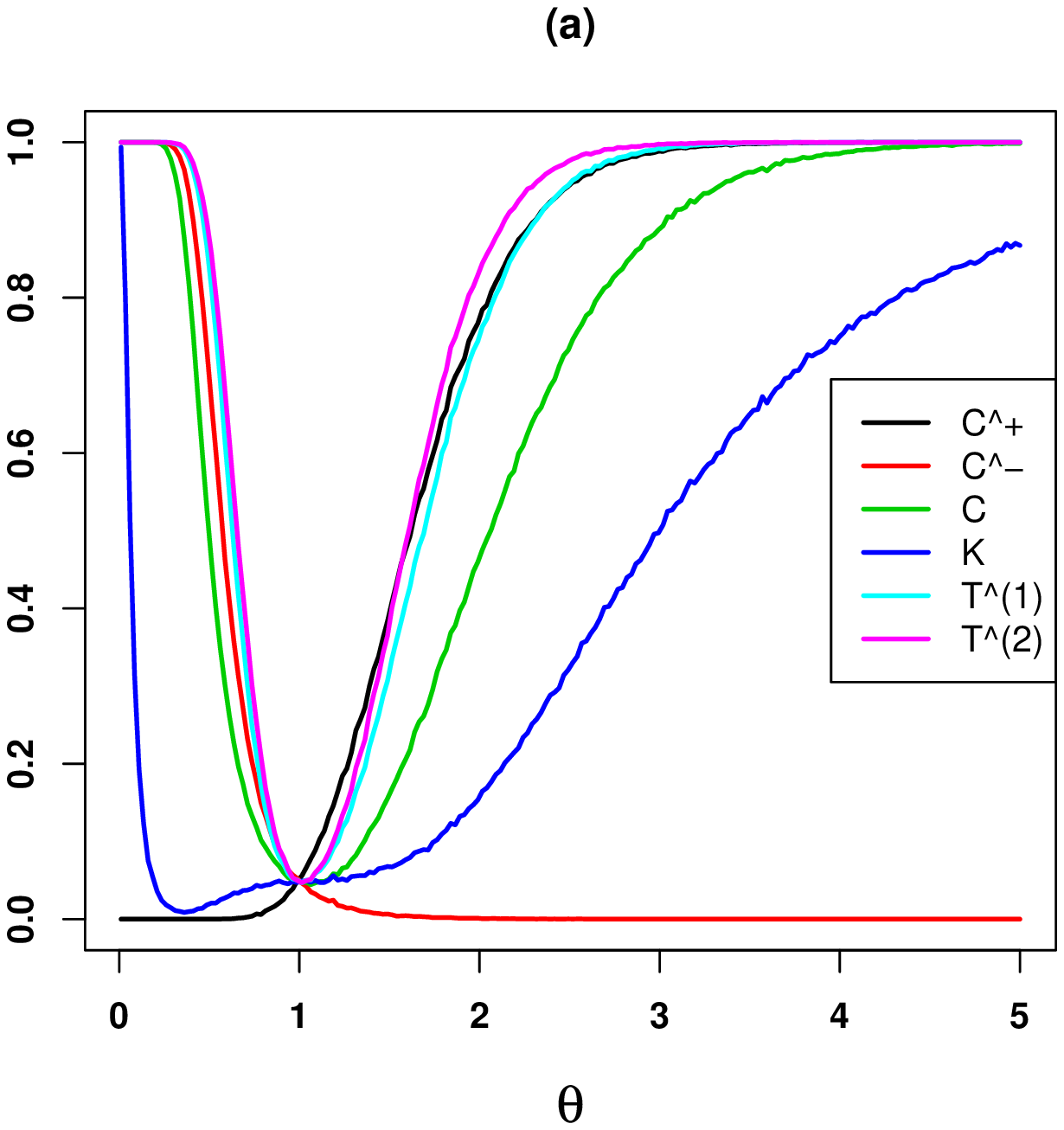} 
\includegraphics[width=4.5cm,height=5cm]{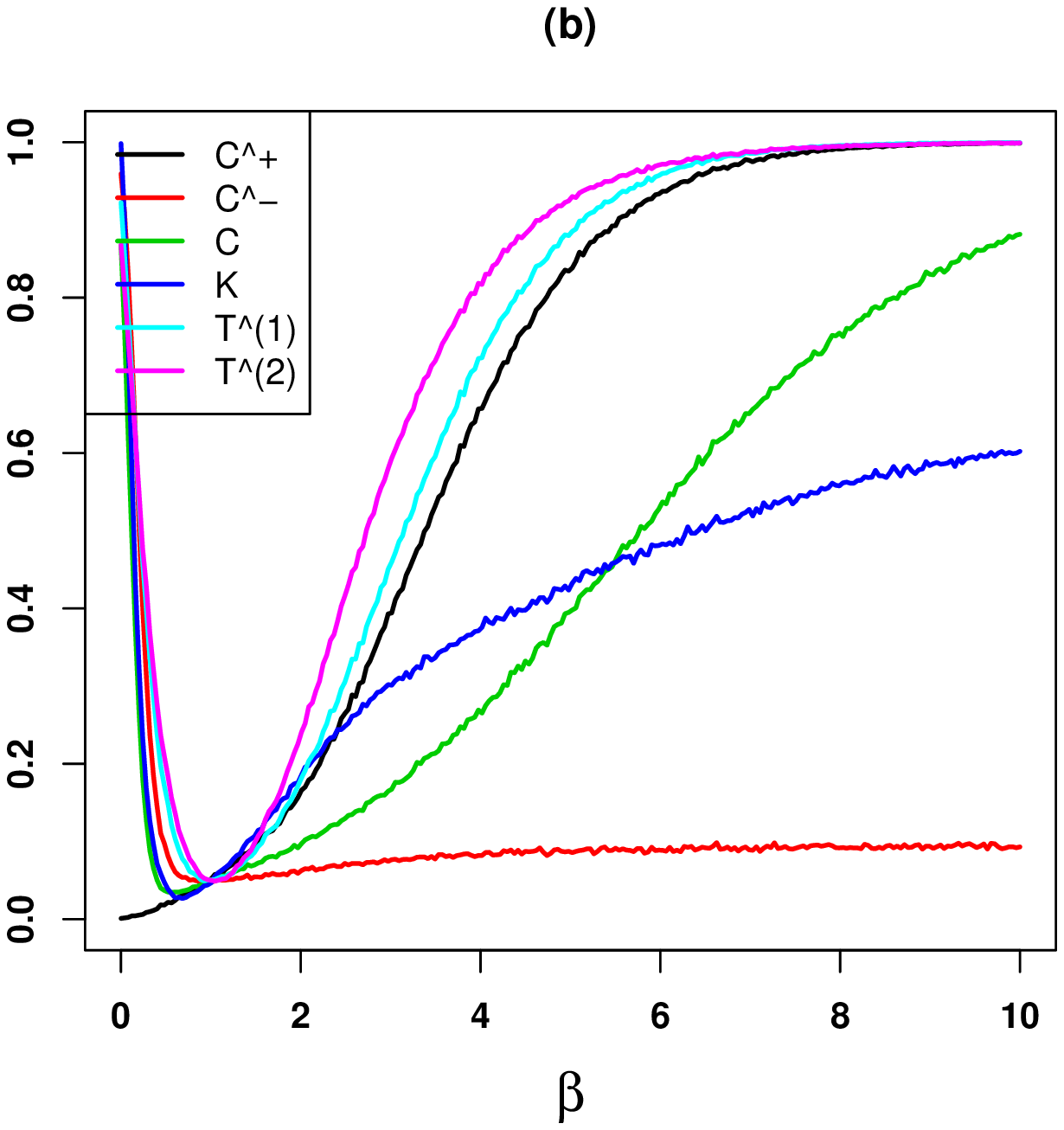}
\includegraphics[width=4.5cm,height=5cm]{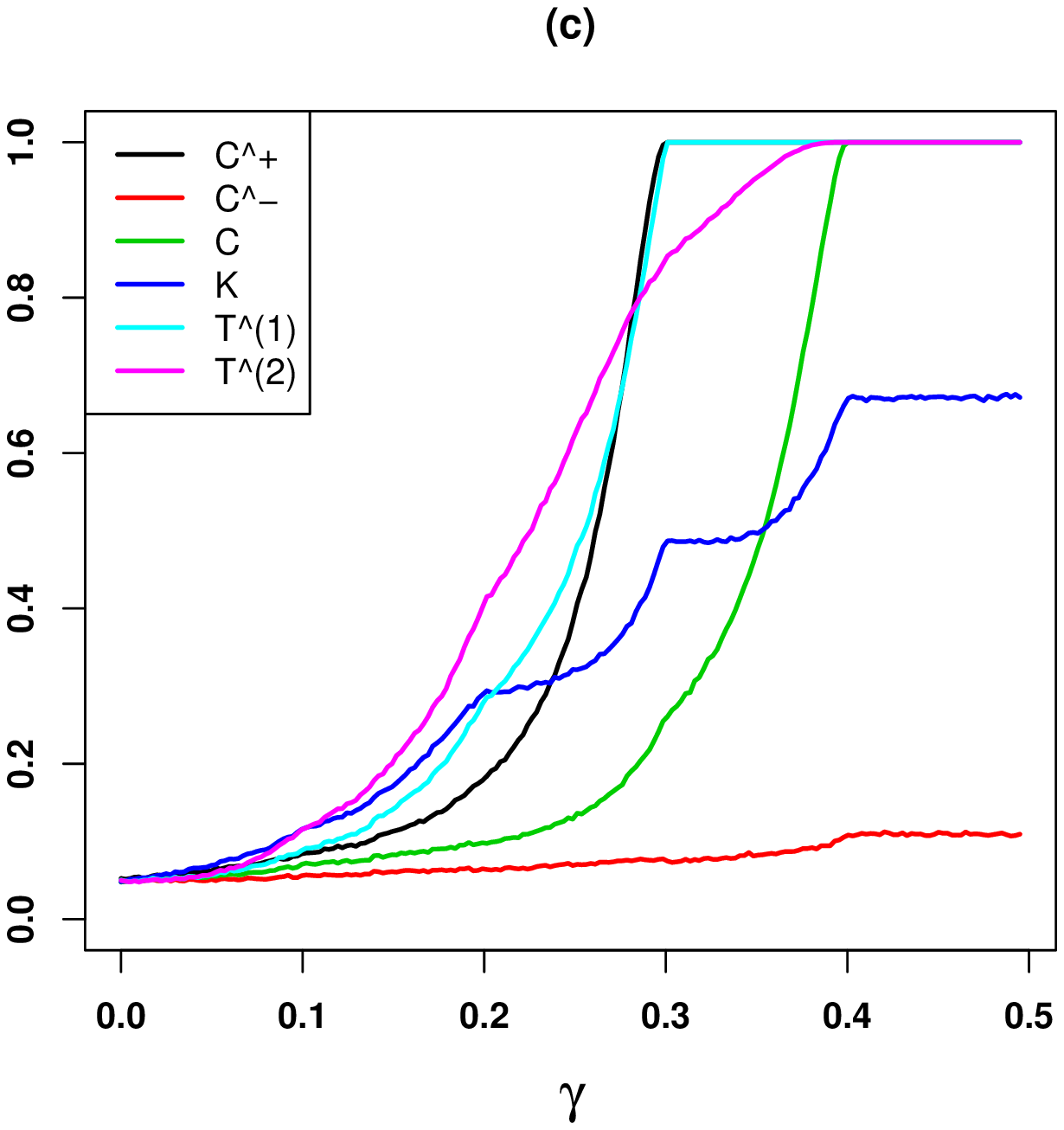}
\caption{\footnotesize{Power function for $\bt=\bt_1$, $\bp=\bp_1$ for families (a), (b) and (c).}}\label{fig1} 
\end{figure}

\begin{figure}[!th] 
\includegraphics[width=4.5cm,height=5cm]{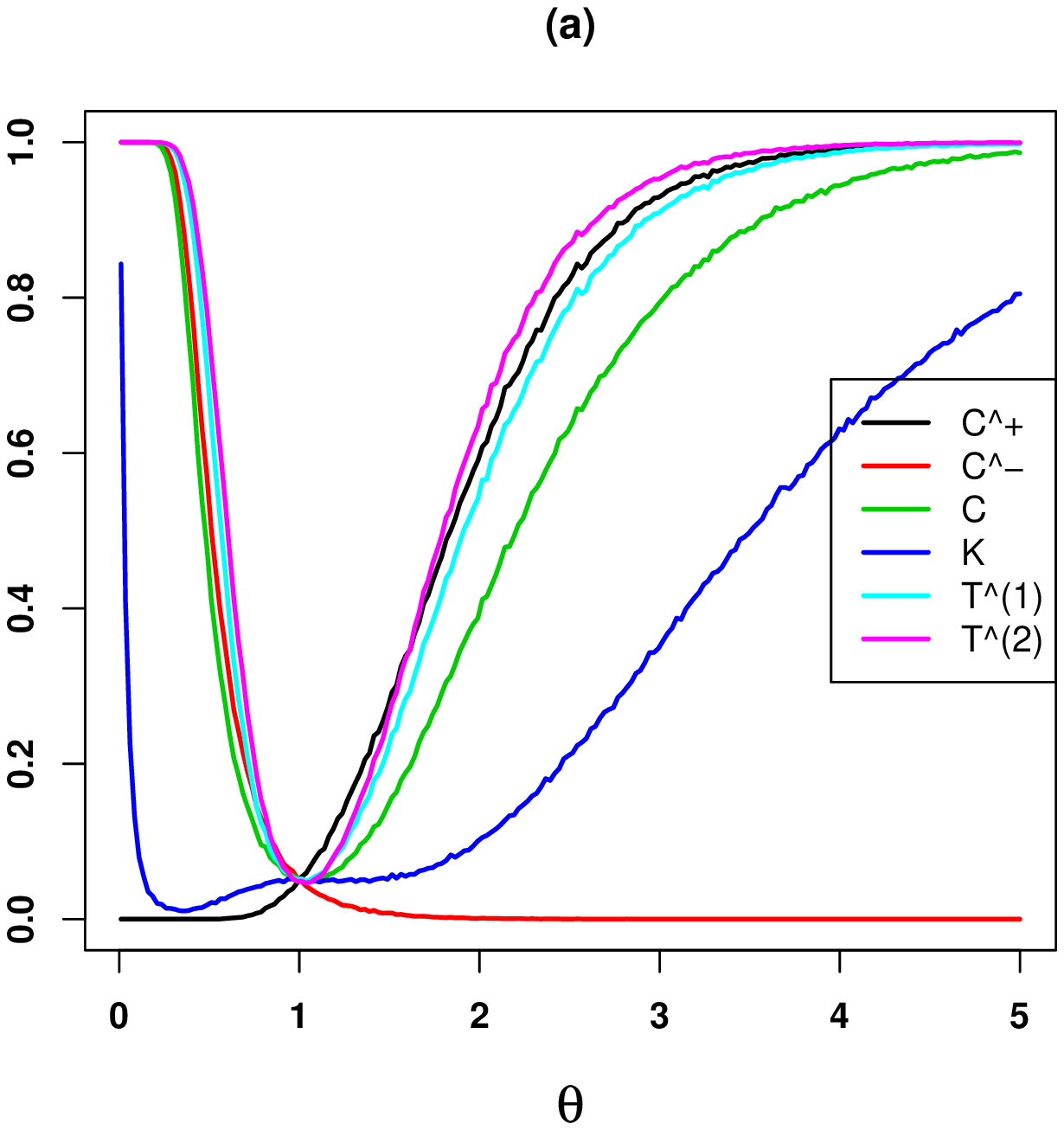} 
\includegraphics[width=4.5cm,height=5cm]{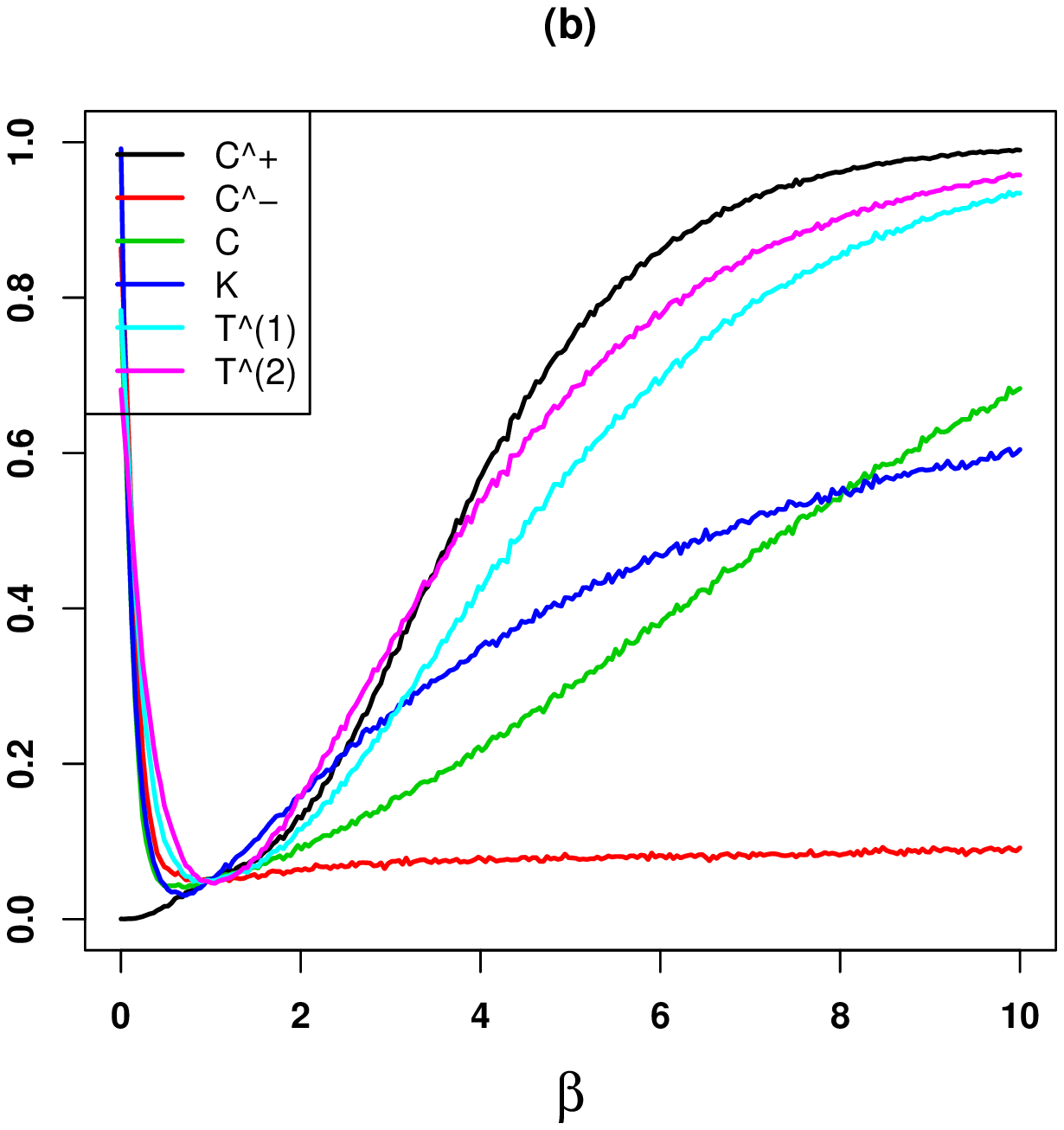}
\includegraphics[width=4.5cm,height=5cm]{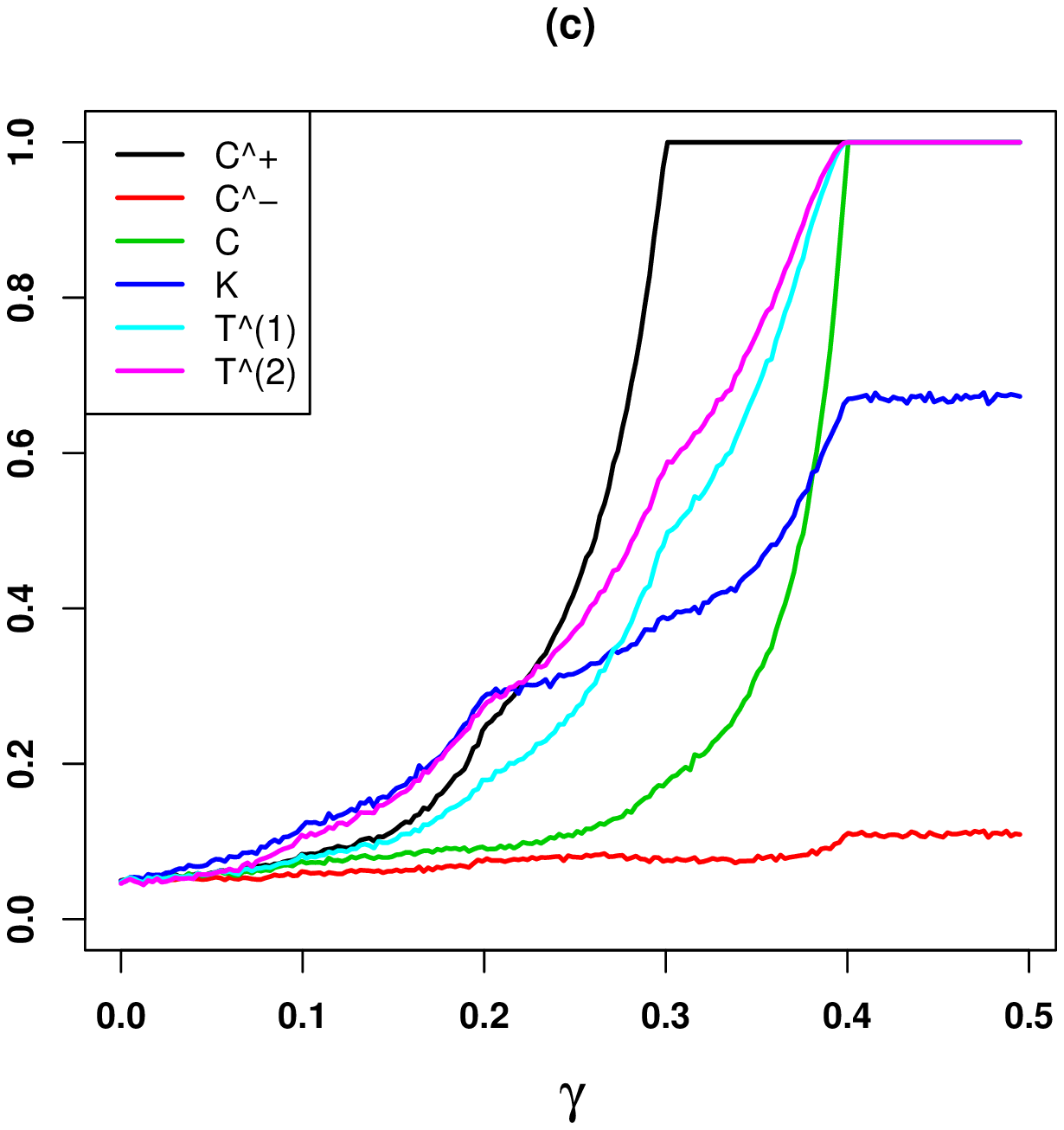}
\caption{\footnotesize{Power function for $\bt=\bt_1$, $\bp=\bp_2$ for families (a), (b) and (c).}}\label{fig2} 
\end{figure}

\begin{figure}[!th] 
\includegraphics[width=4.5cm,height=5cm]{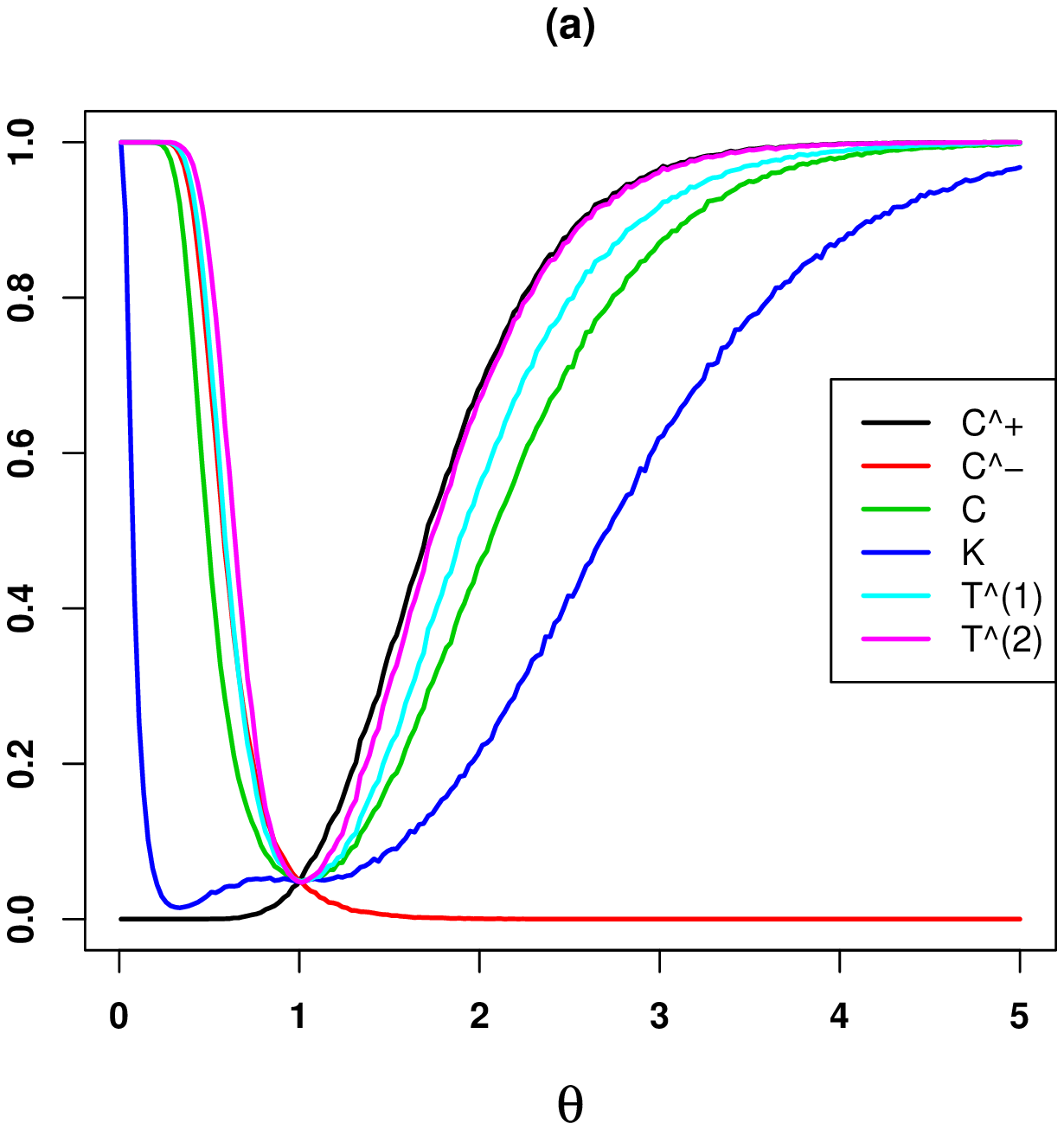} 
\includegraphics[width=4.5cm,height=5cm]{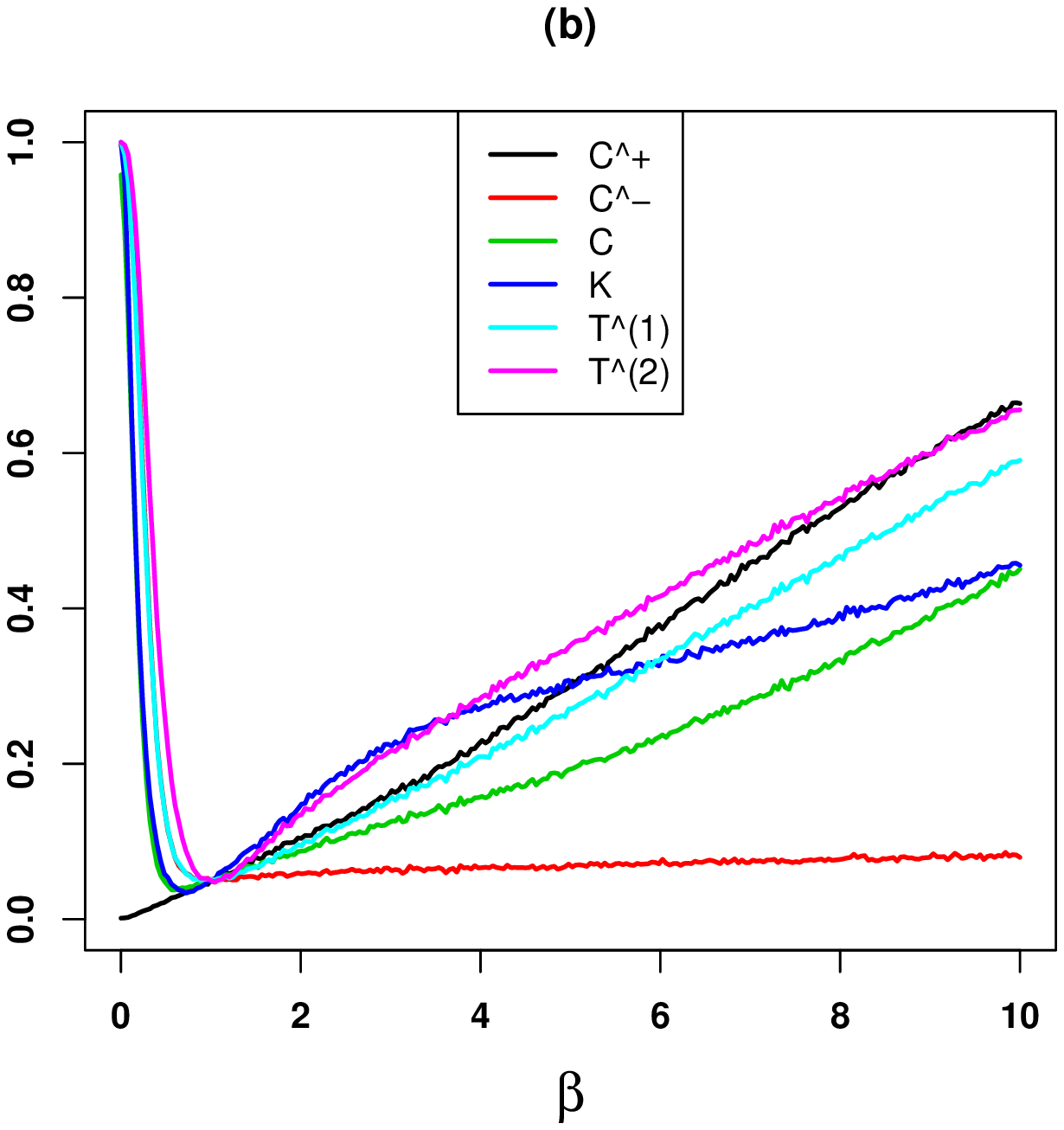}
\includegraphics[width=4.5cm,height=5cm]{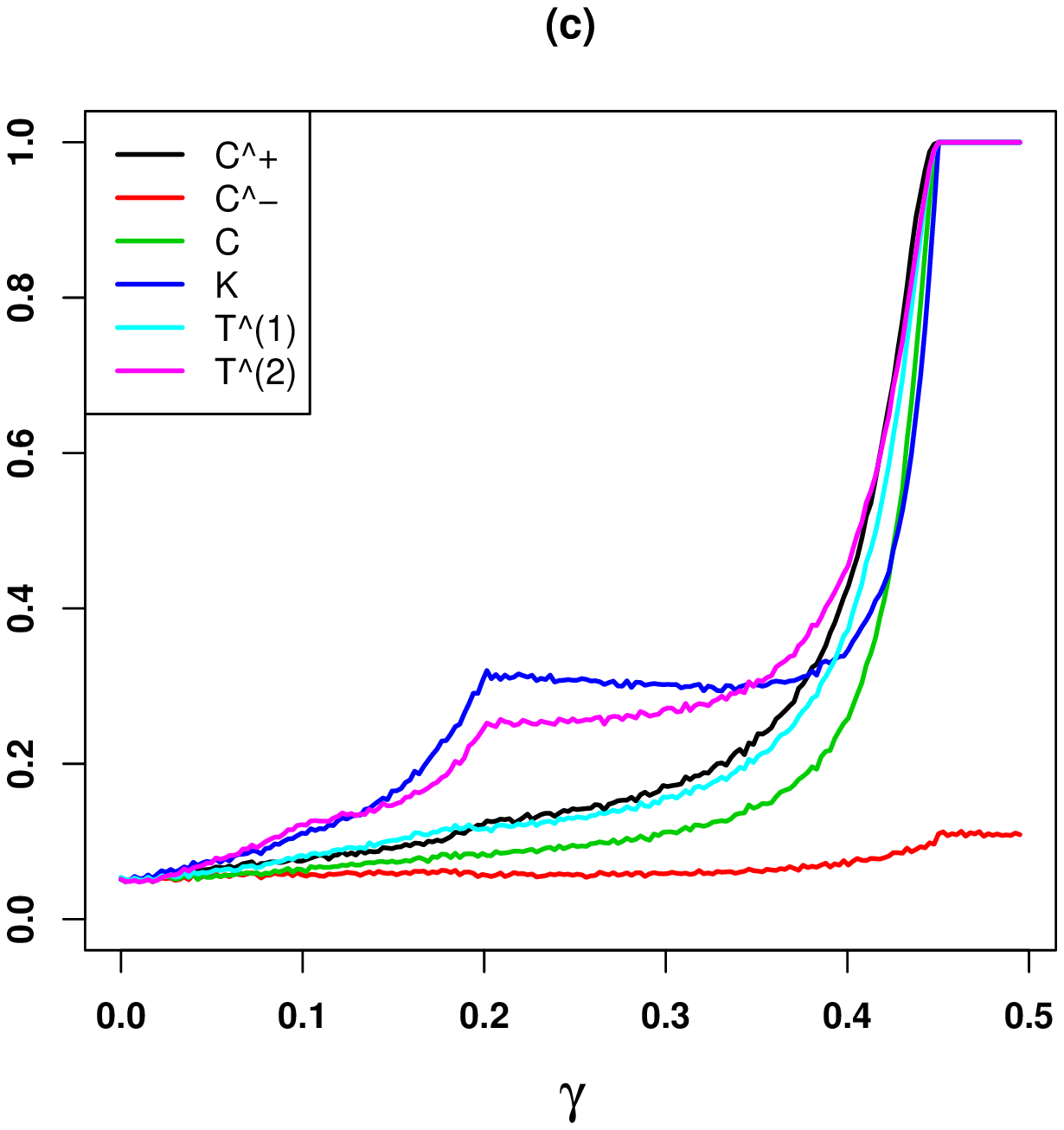}
\caption{\footnotesize{Power function for $\bt=\bt_2$, $\bp=\bp_1$ for families (a), (b) and (c).}}\label{fig3} 
\end{figure}

\begin{figure}[!th] 
\includegraphics[width=4.5cm,height=5cm]{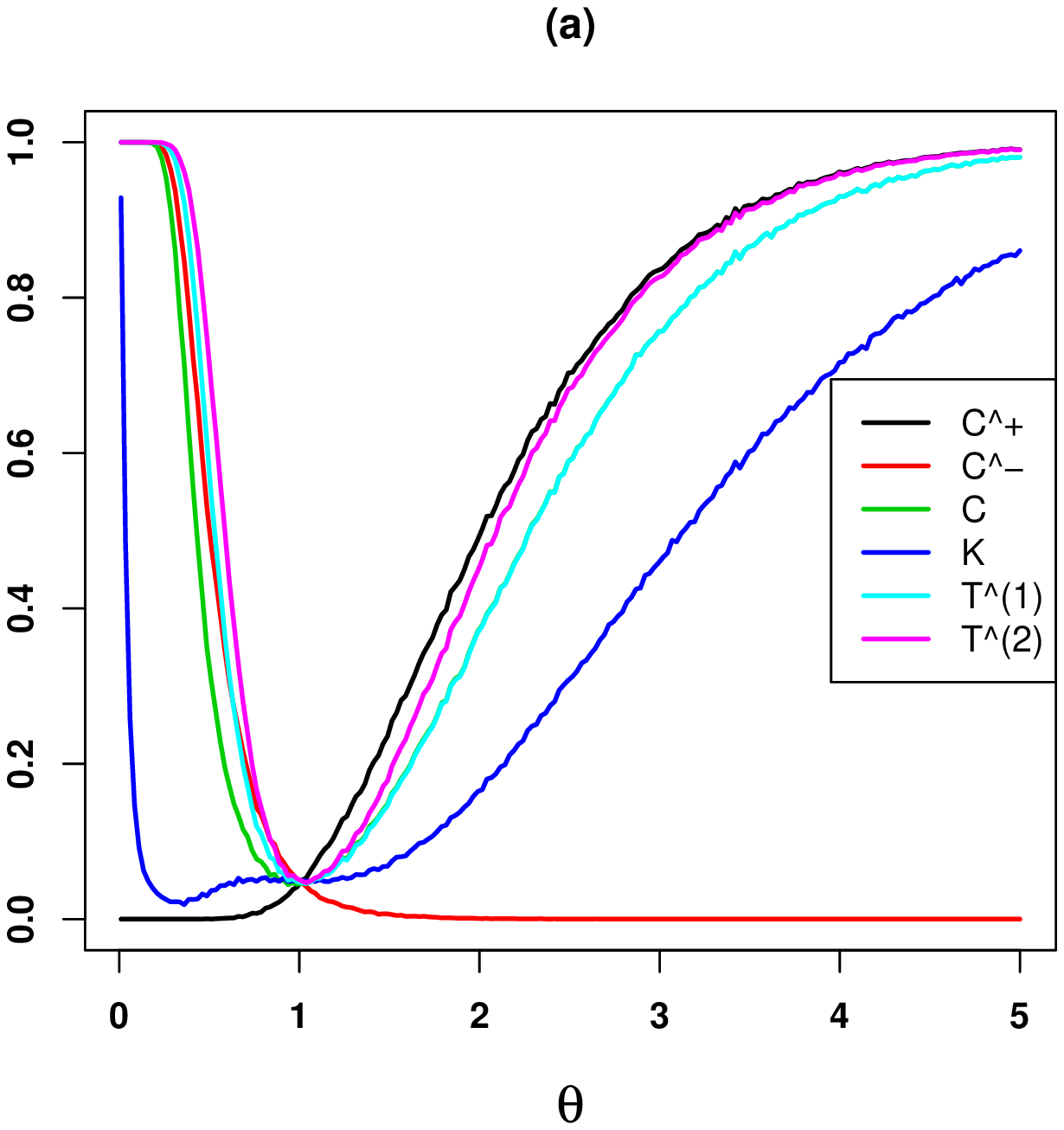} 
\includegraphics[width=4.5cm,height=5cm]{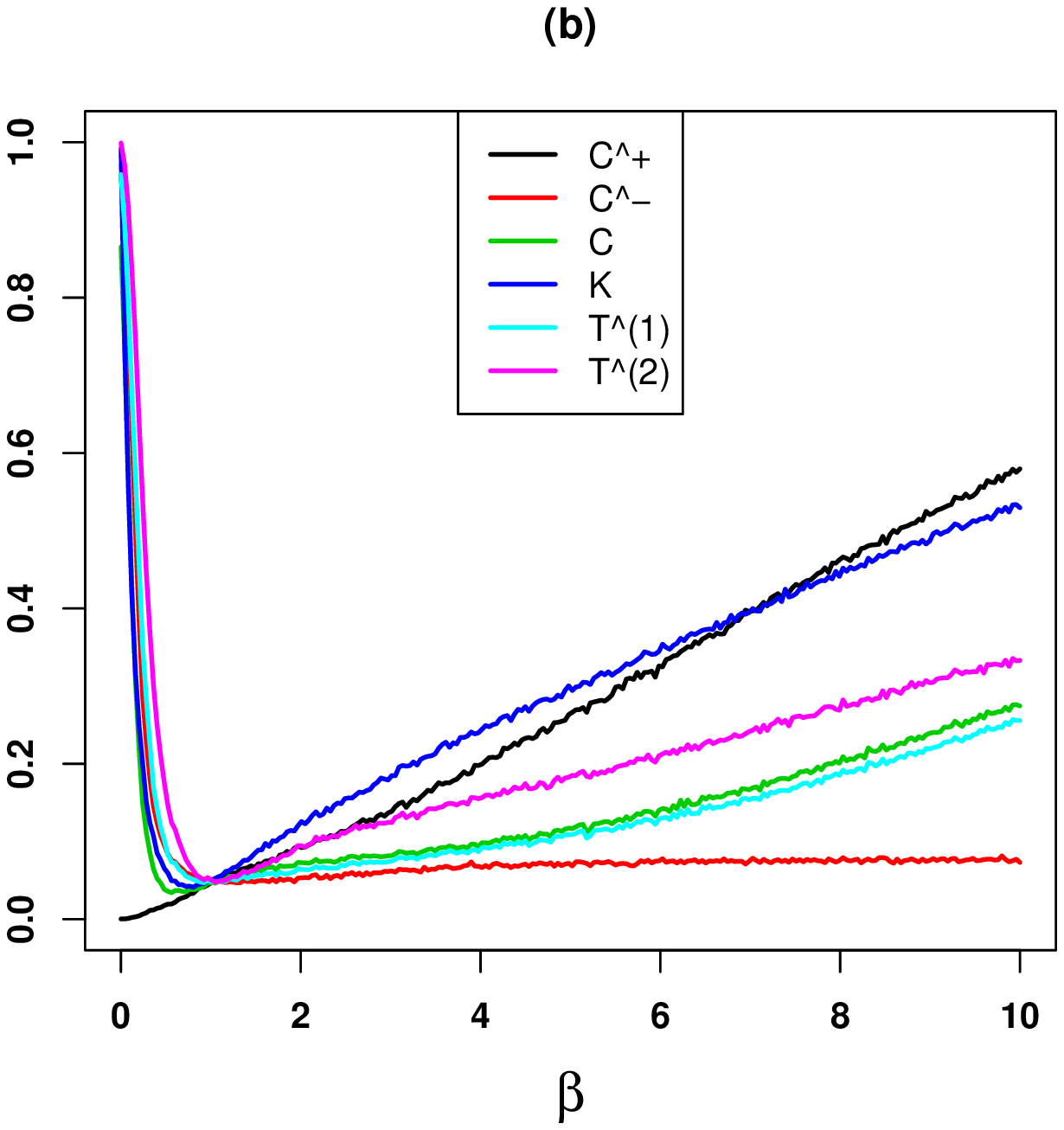}
\includegraphics[width=4.5cm,height=5cm]{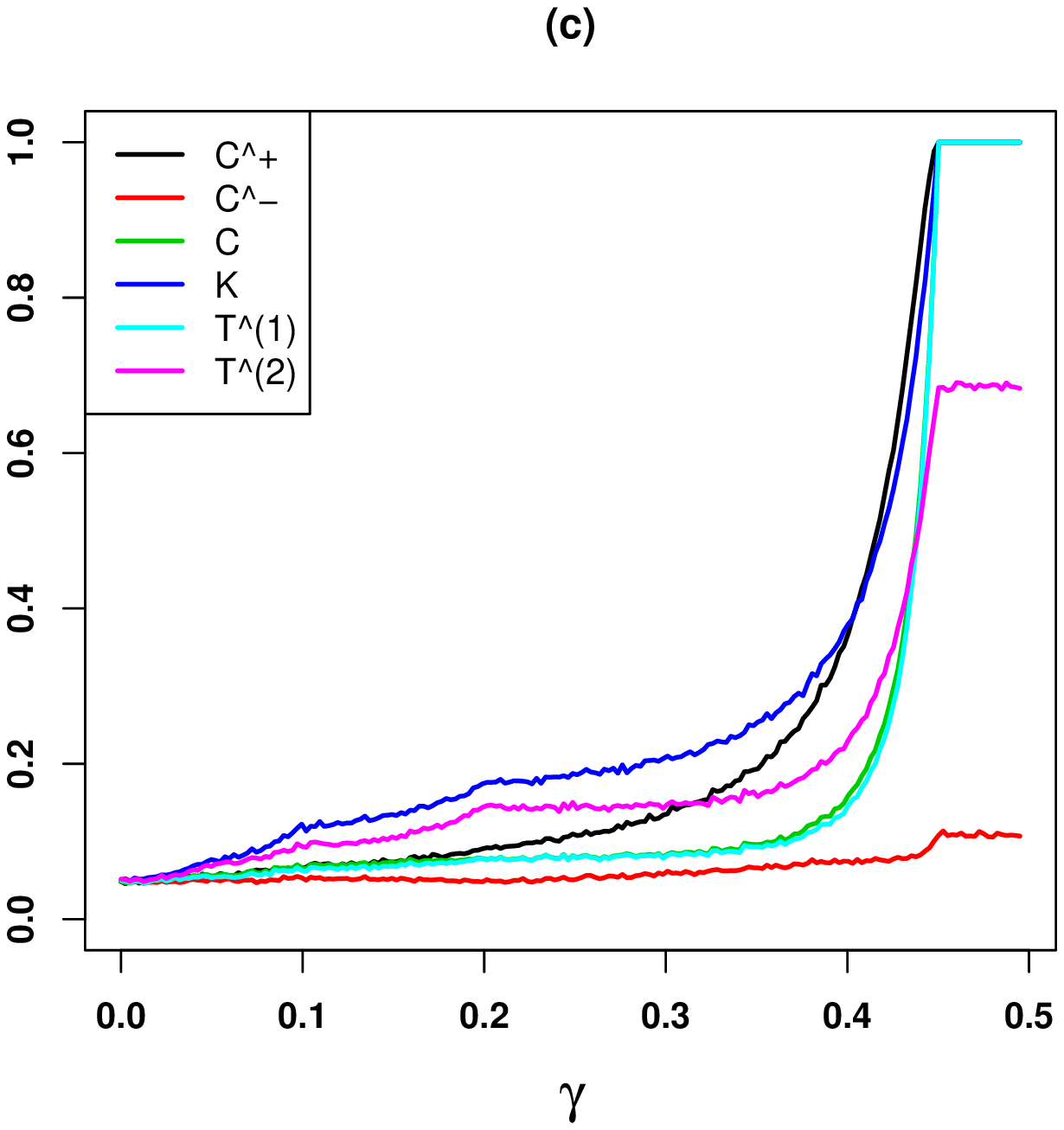}
\caption{\footnotesize{Power functions for $\bt=\bt_2$, $\bp=\bp_2$ for families (a), (b) and (c).}}\label{fig4} 
\end{figure}
\subsection{Discussion}
For family (a), the tests based on $C^+$, $C^-$ ank $K$, and for family (b), the tests based on $C^+$, $C^-$ and $C$ are biased. The tests based on $T^{(1)}$ and $T^{(2)}$ are unbiased. According to the  Figures 1-4, it is clear that, power of the proposed tests depend on the censoring schemes, so we cant find the best test; but it seems that the tests based on $T^{(2)}$, $T^{(1)}$, $C$ and $K$ have good performance, respectively.  

\section{Generalization}
Let $(X_i, R_i, t_i), i=1, \ldots , m$, be progressively Type-I interval censored sample with pre-specified vector $\bp=(p_1, \ldots, p_{m-1}, 1)$ from an unknown distribution function $F(.)$. We are interested in hypothesis testing
\begin{equation}\label{hypothesis2}
\left\{\begin{array}{l}H_0:~F(x)=F_0(x) \\ H_1:~F(x)\neq F_0(x),\end{array}\right.
\end{equation}
where $F_0(.)$ is a continuous and completely specified distribution function. In this case we know that $(X_i, R_i, F_0(t_i)), i=1, \ldots , m$, is a progressively Type-I interval censored sample with pre-specified vector $\bp=(p_1, \ldots, p_{m-1}, 1)$ from $U(0,1)$ distribution. Thus we can test \eqref{hypothesis} by using $(X_i, R_i, F_0(t_i))$, $i=1, \ldots , m$ and $\bp$.

\section{Conclusion}
In this paper, we proposed several statistics for testing uniformity under progressive Type-I interval censoring. We obtained the critical points of these statistics and studied power of the proposed tests  against a representative set of alternatives using simulation. Finally we generalized these methods for  continuous and completely specified distributions.

\end{document}